\documentclass[12pt]{amsart}
\usepackage{amsmath,amssymb,amscd,epsfig,color}
\usepackage{epsfig}
\usepackage{pinlabel}
\newtheorem{theorem}{Theorem}
\newtheorem{corollary}{Corollary}
\newtheorem{degreetheorem}{Degree Theorem}

\def\auts{{$Aut (F_n)$}}

\def\Z{\mathbb Z}
\def\G{\Gamma}

\begin{document}

\title{A Presentation for $Aut(F_n)$}

\begin{author}
{Heather Armstrong, Bradley Forrest, and Karen Vogtmann}
\end{author}

\begin{abstract} We study the action of the group \auts\   of automorphisms of a finitely generated free group on  the degree 2 subcomplex of the spine of Auter space. Hatcher and Vogtmann showed that this subcomplex is simply connected, and we use the method described by K. S. Brown to deduce a new presentation of \auts.  
\end{abstract}
\maketitle

\section{Introduction}

In 1924  Nielsen produced the first finite presentation for the group \auts\ of automorphisms of a finitely-generated free group \cite{nielsen}.  Other presentations have been given by  B. Neumann  \cite{neumann} and J. McCool  \cite{mccool}.   A  very natural presentation for the index two subgroup $SAut(F_n)$ was given by Gersten in \cite{gersten}.

Nielsen, McCool and Gersten used infinite-order generators.  Neumann used only finite-order generators of order at most $n$, but his relations are very complicated.  P. Zucca showed that $Aut(F_n)$ can be generated by three involutions, two of which commute, but did not give a complete presentation \cite{zucca}. 

In this paper we produce a new presentation for $Aut(F_n)$ which has several interesting features.    The generators are involutions and the number of relations is fairly small.  The form of the presentation for $n\geq 4$ depends only on the size of a signed symmetric subgroup.  
 
The presentation is found by considering the action of $Aut(F_n)$ on a subcomplex of   the {\it spine of Auter space}.  This spine is a contractible simplicial complex on which $Aut(F_n)$ acts with finite stabilizers and finite quotient.  A vertex  of the spine corresponds to a basepointed graph $\G$ together with an isomorphism $F_n\to \pi_1(\G)$.   In \cite{hatcher-vogtmann}  Hatcher and Vogtmann defined a sequence of nested invariant subcomplexes $K_r$ of this spine, with the property that the  $r$-th complex $K_r$ is  $(r-1)$-connected.  In particular, $K_2$ is simply-connected, and we use  the method described by K. S. Brown in \cite{brown} to produce our finite presentation using the action of $Aut(F_n)$ on $K_2$.

In order to describe the presentation, we 
fix generators $a_1,\ldots,a_n$ for the free group $F_n$  and let $W_n$  be the subgroup of $Aut(F_n)$ which permutes and inverts the generators.   We let $\tau_i$ denote the element of $W_n$ which inverts $a_i$, and $\sigma_{ij}$ the element which interchanges $a_i$ and $a_{j}$: $$
{\tau_i\colon\begin{cases}  a_i\mapsto a_i^{-1}&\cr
            a_j\mapsto a_j & j\neq i
            \end{cases}}
\qquad\qquad{\sigma_{ij}\colon\begin{cases}  a_i\mapsto a_j\cr
            a_j\mapsto a_i\cr
            a_k\mapsto a_k& k\neq i,j.\end{cases}}$$
There are many possible presentations of $W_n$.  For instance,   
$W_n$ is generated by $\tau_1$   and by transpositions $s_i=\sigma_{i,i+1}$ for $1\leq i\leq n-1$,             
subject to relations 
 \begin{eqnarray*}
s_{i}^2=1 & 1\leq i\leq n-1\\
(s_i s_j)^2= 1  &  j\neq i\pm 1\\
 (s_{i-1}s_{i})^3=1  &2\leq i\leq i-1\\
 \tau_1^2=1\\
 (\tau_1 s_{1})^4 = 1\\
 \left(\tau_1 s_{i}\right)^2=1 & 2\leq i\leq i-1.
\end{eqnarray*}
Generators for $Aut(F_n)$ will consist of generators for $W_n$ plus the following involution:
\begin{equation*}
\eta\colon \begin{cases}a_1\mapsto{a_2^{-1}}a_1\\
a_2\mapsto a_2^{-1} \\
a_k\mapsto a_k  & k>2.
\end{cases}
\end{equation*}
The presentation we obtain is the following:
\begin{theorem}\label{presentation}  For $n \geq 4$, $Aut(F_n)$ is generated by $W_n$ and $\eta$, subject to the following relations:
\begin{enumerate}
\item$\eta^2=1$
\item$ (\sigma_{12}\eta)^3=1$
\item$ (\eta \tau_{i})^2= 1$ for $ i>2$
\item$(\eta \sigma_{ij})^2= 1$ for $ i,j>2$
\item$((\eta\tau_1)^2\tau_2)^2=1$
\item $(\eta\sigma_{13}\tau_2\eta\sigma_{12})^4=1$
\item $ \sigma_{12}\eta\sigma_{13}\tau_2\eta\sigma_{12}(\sigma_{23}\eta\sigma_{13}\tau_2\eta)^2= 1$
\item  $(\sigma_{14}\sigma_{23}\eta)^4=1$ 
\item relations in $W_n$.
\end{enumerate}
\end{theorem}
The presentation we obtain for $n=3$ differs only in that every relation involving indices greater than 3 is missing:  
\begin{corollary} The group $Aut(F_3)$ is generated by $W_3$ and $\eta$, subject to the following relations:
\begin{enumerate}
\item$\eta^2=1$
\item$ (\sigma_{12}\eta)^3=1$
\item$ (\eta \tau_{3})^2= 1$  
\item$((\eta\tau_1)^2\tau_2)^2=1$
\item $(\eta\sigma_{13}\tau_2\eta\sigma_{12})^4=1$
\item $  \sigma_{12}\eta\sigma_{13}\tau_2\eta\sigma_{12}(\sigma_{23}\eta\sigma_{13}\tau_2\eta)^2= 1$
\item relations in $W_3$.
\end{enumerate}
\end{corollary}
 
For $n=2$ we get:
\begin{corollary} The group $Aut(F_2)$ is generated by $\tau_1,\tau_2,\sigma_{12}$ and $\eta$, subject to the following relations:
\begin{enumerate}
\item$\eta^2=1$
\item$ (\sigma_{12}\eta)^3=1$
\item$((\eta\tau_1)^2\tau_2)^2=1$
\item $\sigma_{12}^2=1$
\item $\tau_1^2=1$
\item $(\tau_1\sigma_{12})^4=1$
\item $\tau_2=\sigma_{12}\tau_1\sigma_{12}$.
\end{enumerate}
\end{corollary}

\section{Brown's theorem}\label{brownstheorem}

To find  our presentation, we use the method described by K. S. Brown in \cite{brown}.  This method applies whenever a group G acts on a simply-connected CW-complex by permuting cells, but the desciption is simpler if the complex is simplicial and the action does not invert edges.  Since this is the case for us, we describe this simpler version.    We remark that a presentation of the fundamental group of a complex of groups, whether or not it arises from the action of a group on a complex, can be found in \cite{BH}, Chapter III.$\mathcal C$.  

Let $G$ be a group, and $X$ a non-empty simply-connected simplicial complex on which $G$ acts without inverting any edge.  Let $\mathcal{V}, \mathcal{E}$ and $\mathcal{F}$  be sets of representatives of vertex-orbits, edge-orbits, and 2-simplex-orbits, respectively, under this action.   The group $G$ is generated by the stabilizers $G_v$ of vertices in $\mathcal{V}$ together with a generator for each edge $e\in\mathcal{E}$.    There is a relation for each element of $\mathcal{F}$.  Other relations come from loops in the 1-skeleton of the quotient $X/G$.  In order to write down a presentation explicitly, we choose the sets $\mathcal{V}$, $\mathcal{E}$ and $\mathcal{F}$ quite carefully, as follows.

The 1-skeleton of the quotient $X/G$ is a graph.  Choose a maximal tree in this graph  and lift it to a tree $\mathcal{T}$ in $X$.  The vertices of $\mathcal{T}$ form $\mathcal{V}$,  our set of vertex-orbit representatives for the action of $G$ on   $X$.  Since the edges of $\mathcal{T}$ are not a complete set of edge-orbit representatives, we complete the set $\mathcal{E}$ by including for each missing orbit a choice of representative which is connected to $\mathcal{T}$. Finally, for the set $\mathcal{F}$, we choose representatives for the 2-simplices so that they also have at least one vertex in $\mathcal{T}$.
 
We obtain a presentation for $G$ as follows:

\noindent{\bf Generators}.    The group $G$ is generated by the stabilizers $G_v$ for $v\in\mathcal{V}$ together with a generator $t_e$ for each $e\in \mathcal{E}$.

\noindent{\bf Relations}.  There are four types of relations:  tree relations, edge relations, face relations and stabilizer relations.  The {\it tree relations} are:
\begin{enumerate}
\item $t_e=1$ if $e\in \mathcal{T}.$
\end{enumerate}
There are {\it edge relations} for each edge $e\in \mathcal{E}$ which identify the two different copies of $G_e$, the stabilizer of e, which can be found in the stabilizers of the endpoints of $e$.  To make this explicit, we orient each edge $e\in \mathcal{E}$ so that the initial vertex $o(e)$ lies in $\mathcal{T}$, and let $i_e: G_e \to G_{o(e)}$ denote the inclusion map.  There is also an inclusion $G_e \to G_{t(e)}$, where $t(e)$ is the terminal vertex.  Note that when $t(e)$ is not in $\mathcal{T}$, $G_{t(e)}$ is not in our generating set.  To encode the information of this inclusion map in terms of our generating set we must do the following.  Since $t(e)$ is equivalent to some vertex $w(e)$ in $\mathcal{T}$, we choose $g_e \in G$ with $g_ew(e)=t(e)$ (if $t(e)\in \mathcal{T}$, we choose $g_e=1$).   Conjugation by $g_e$ is an isomorphism from $G_{t(e)}$ to $G_{w(e)}$, so we set $c_e\colon  G_e\to G_{w(e)}$ to be the inclusion $G_e\to G_{t(e)}$ followed by conjugation by $g_e$.   Equating the two images of $G_e$ gives us the edge relations, which are then:  

\begin{enumerate}\addtocounter{enumi}{1}
\item For $x\in G_e$,\,\,  $t_e i_e(x)t_e^{-1}=c_e(x)$.
\end{enumerate}
There is a {\it face relation} for each $2$-simplex $\Delta\in\mathcal{F}$.  To describe this, we use the notation established in the previous paragraph.   

We digress for a moment to consider an arbitrary oriented edge  $e'$ of $X$ with $o(e') \in \mathcal{V}$.  This edge is equivalent to some edge $e\in\mathcal{E}$.  If the orientations on $e'$ and $e$ agree, then  
  $e'=he$ for some $h\in G_{o(e')}$, and  $t(e')= hg_ew(e)$.
If the orientations do not agree, then $e'=h g_e^{-1}e$ for some $h\in G_{o(e')}$, and $t(e')= h g_e^{-1}o(e)$.
The element  $h$ is unique modulo the stabilizer of $e'$.  

Now let $e'_1e'_2e'_3$ be an oriented edge-path starting in $\mathcal{T}$ and going around the boundary of $\Delta$.  Since $e'_1$ originates in $\mathcal{T}$, we can associate to it elements $h_1\in G_{o(e_1')}$ and $g_1=h_1g_{e_1}^{\pm1}$  as described above. 
Then $e'_2$ 
originates in $g_1\mathcal{T}$, so $g_1^{-1}e'_2$ originates in $\mathcal{T}$, and we can find $h_2$ and    $g_2=h_2g_{e_2}^{\pm1}$  for $g_1^{-1}e'_2$.  Now
$e'_3$ originates in $g_1 g_2 \mathcal{T}$ so we can find $h_3$ and   $g_3=h_3g_{e_3}^{\pm1}$   associated to  $g_2^{-1}g_1^{-1}e'_3$.  Set $g_\Delta=g_1g_2g_3$, and note that $g_\Delta$ is in the stabilizer of the vertex $o(e_1')$, so that the following is a relation among our generators:
\begin{enumerate}\addtocounter{enumi}{2}
\item For each $\Delta\in\mathcal{F}, \,\,  h_1t_{e_1}^{\pm1}h_2t_{e_2}^{\pm1}h_3t_{e_3}^{\pm1}=g_\Delta$.
\end{enumerate}
Here the sign on $t_{e_i}$ is equal to the sign on $g_{e_i}$ in the expression for $g_i$.

Finally, a {\it stabilizer relation} is a relation among the generators of a vertex stabilizer $G_v$.

\medskip
\begin{theorem}\label{brown}
 {{\bf (Brown)}} Let $X$ be a simply-connected simplicial complex with a simplicial action by the group $G$ which does not invert edges.   Then $G$ is generated by the stabilizers $G_v$ ($v\in \mathcal{V}$) and symbols $t_e$ ($e\in \mathcal{E}$) subject to all tree, edge, face and stabilizer relations as described above.\\
\end{theorem}

\section{Computations} 
\medskip
\subsection{The Degree 2 complex}
We will apply Theorem ~\ref{brown}  to a certain subcomplex of the spine of Auter space.  The spine of Auter space is a contractible simplicial complex on which $Aut(F_n)$ acts with finite stabilizers and finite quotient.    For full details on the construction of Auter space, we refer to \cite{hatcher-vogtmann}.

\begin{figure}[ht!]
\labellist
\small\hair 2pt  
\pinlabel {$a_1$} [r] at -1 75
\pinlabel {$a_3$} [r] at 1 14
\pinlabel {$a_2$} [l] at 65 75
\pinlabel {$a_4$} [l] at 60  14
\pinlabel {$=$}  at 100 54
\pinlabel {$a_2^{-1}a_1$} [r] at 175 75
 \pinlabel {$a_3$} [r] at 180 14
 \pinlabel {$a_2$} [l] at 214 75
 \pinlabel {$a_4$} [l] at 235  14
\endlabellist
\centering
\includegraphics[scale=0.75]{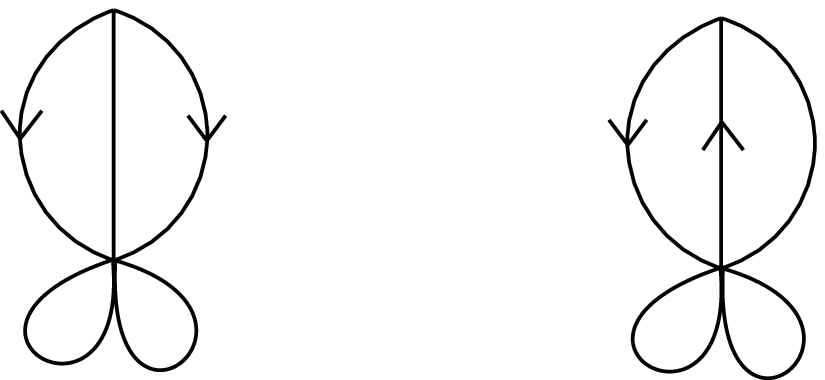}
\caption{Two labeled graphs representing the same vertex in the spine of Auter space}
\label{point}
\end{figure}

A vertex in the spine of Auter space is a connected, basepointed graph $\Gamma$ together with an isomorphism $g\colon\pi_1(\Gamma)\to F_n$, called a {\it marking}.   (Note: often in the literature the marking goes in the other direction).  We require all vertices of $\Gamma$ to have valence at least three, and we also assume that $\Gamma$ has no separating edges.  One can describe the marking $g$ by labeling certain edges of $\Gamma$ as follows.   Choose a maximal tree in $\Gamma$. The edges not in this maximal tree form a natural basis for the fundamental group $\pi_1(\Gamma)$.  Orient each of these, and label them by their images in $F_n$.  This description depends on the choice of maximal tree; for instance the labeled graphs in Figure~\ref{point} represent the same vertex.  
 
Two marked graphs span an edge in the spine if one can be obtained from the other by collapsing a set of edges (this is called a {\it forest collapse}).  A set of $k+1$ vertices spans a $k$-simplex if each pair of vertices spans an edge.

\begin{figure}[ht!]
\labellist
\small\hair 2pt 
\pinlabel{$\eta\cdot$} [r] at -30 54 
\pinlabel {$a_1$} [r] at -1 75
\pinlabel {$a_3$} [r] at 1 14
\pinlabel {$a_2$} [l] at 65 75
\pinlabel {$a_4$} [l] at 60  14
\pinlabel {$=$}  at 100 54
\pinlabel {$a_2^{-1}a_1$} [r] at 163 75
 \pinlabel {$a_3$} [r] at 170 14
 \pinlabel {$a_2^{-1}$} [l] at 233 75
 \pinlabel {$a_4$} [l] at 225  14
\endlabellist
\centering
\includegraphics[scale=0.75]{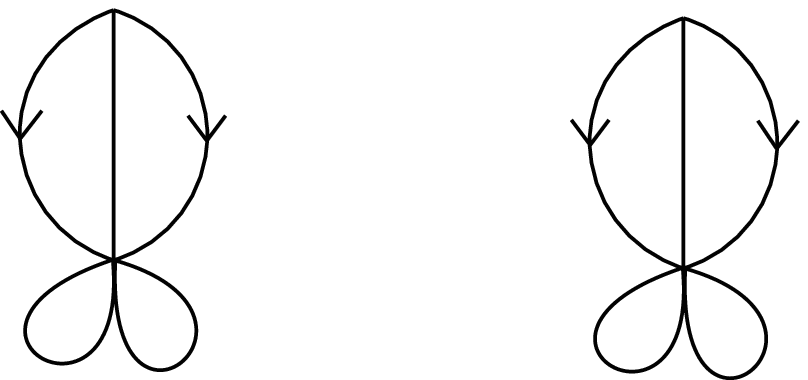}
\caption{Action of $\eta$ on a Nielsen graph.  }
\label{action}
\end{figure}

The group $Aut(F_n)$ acts on Auter space on the left by $\alpha\cdot(g,\Gamma)=(\alpha\circ g,\Gamma)$.  This is represented on a labeled graph by applying $\alpha$ to the edge-labels. Figure~\ref{action} shows the results of applying $\eta$ to the graph from Figure~\ref{point}.  Note that this is the same marked graph, so that $\eta$ fixes this vertex of the spine. 

The {\it degree} of a graph is defined to be  $2n$ minus the valence of the basepoint. The only graph of degree 0 is a rose, and the only graph of degree 1 is the graph underlying the marked graph in Figure~\ref{point}.  There are five different graphs of degree 2.  

 A forest collapse cannot increase degree, so the vertices of degree at most $i$ span a subcomplex $K_i$ of the spine.  Hatcher and Vogtmann proved that the subcomplexes $K_i$ act like ``skeleta" for the spine of Auter space: 

\noindent\begin{degreetheorem}{{\bf   \cite{hatcher-vogtmann}}}  $K_i$ is $i$-dimensional and $(i-1)$-connected.  
\end{degreetheorem}

In particular, the subcomplex $K_2$ spanned by graphs of degree at most 2  is a simply-connected 2-complex.

\subsection{Quotient}

\begin{figure}[ht!]
\labellist
\small\hair 2pt 
\pinlabel{$\not\in \mathcal{E}$} [l] at 60 0
\pinlabel{$\in \mathcal{E\backslash T}$} [l] at 60 14
\pinlabel{$\in \mathcal{T}$} [l] at 60 28
\pinlabel{${\bf v}_2$}   at 33 73
\pinlabel{${\bf v}_5$}   at 330 73
\pinlabel{${\bf v}_0$}   at 200 73
\pinlabel{${\bf v}_3$}   at 43 203
\pinlabel{${\bf v}_6$}   at 318 203
\pinlabel{${\bf v}_1$}   at 168 203
\pinlabel{${\bf v}_4$}   at 80 285
\pinlabel{${\bf v}_6$}   at 288 285
\pinlabel{${\bf v}_8$}   at 233 287
\pinlabel{$\Delta_{104}$}  at 150 170
\pinlabel{$\Delta_{107}$}   at 220 170
\scriptsize
\pinlabel{$a_1$} at 15 90
\pinlabel{$a_2$} at 30 53
\pinlabel{$a_3$} at 23 20
\pinlabel{$a_4$} at 5 20
\pinlabel{$a_1$} at 150 55
\pinlabel{$a_2$} at 205 55
\pinlabel{$a_3$} at 150 20
\pinlabel{$a_4$} at 205 20
\pinlabel{$a_1$} at 335 55
\pinlabel{$a_2$} at 360 55
\pinlabel{$a_3$} at 385 55
\pinlabel{$a_4$} at 360 22
\pinlabel{$a_1$} at 23 273
\pinlabel{$a_2$} at 35 240
\pinlabel{$a_3$} at 10 210
\pinlabel{$a_4$} at 37  220
\pinlabel{$a_1$} at 320 240
\pinlabel{$a_2$} at 354 240
\pinlabel{$a_3$} at 389 240
\pinlabel{$a_4$} at 357  205
\pinlabel{$a_1$} at 195 300
\pinlabel{$a_2$} at 218 325
\pinlabel{$a_3$} at 235 325
\pinlabel{$a_4$} at 260 300
\pinlabel{$a_1$} at 68 325
\pinlabel{$a_2$} at 53 305
\pinlabel{$a_3$} at 45 267
\pinlabel{$a_4$} at 60 267
 \pinlabel{$a_1$} at 295 305
\pinlabel{$a_2$} at 330 325
\pinlabel{$a_3$} at 368 305
\pinlabel{$a_4$} at 332 272
 \pinlabel{$a_1$} at 148 265
\pinlabel{$a_2$} at 192 265
\pinlabel{$a_3$} at 164 223
\pinlabel{$a_4$} at 179 223
\endlabellist
\centering
\includegraphics[scale=0.9]{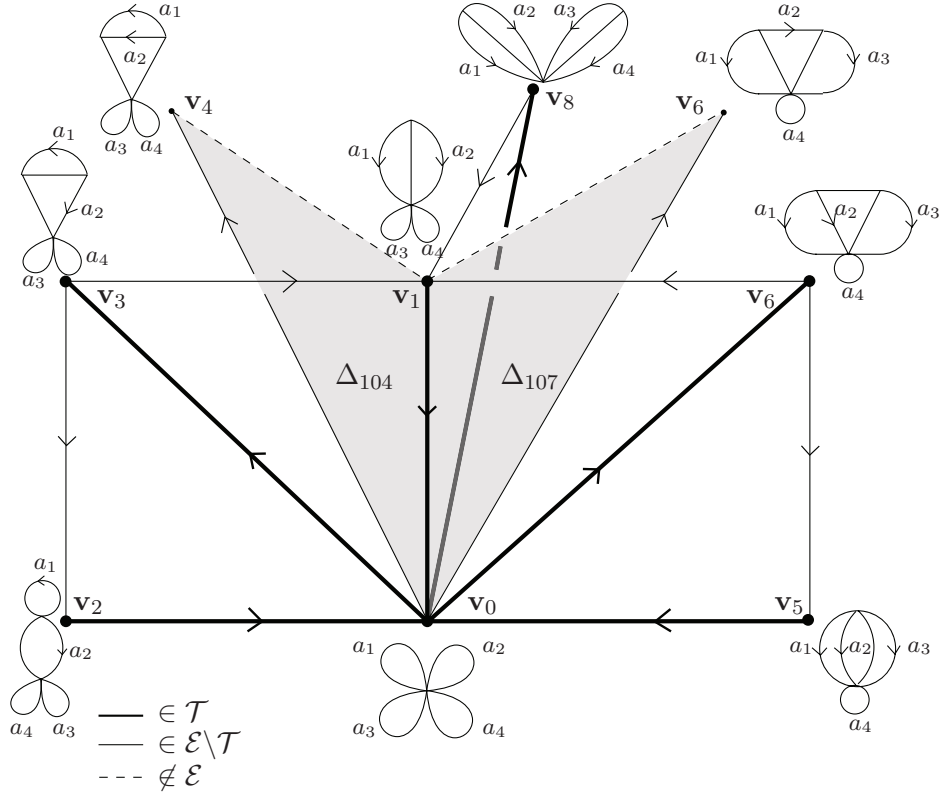}
\caption{A lift of the quotient $K_2/Aut(F_4)$.  For $n>4$, add loops to the basepoint.}
\label{qt}
\end{figure}

The quotient of $K_2$ by the action of $Aut(F_n)$ was computed in \cite{hatcher-vogtmann}.   For $n\geq 4$, this quotient  has seven vertices, thirteen edges and seven triangles.  Figure~\ref{qt} shows a lift of these simplices to $K_2$ for $n=4.$ For $n>4$, the picture is the same except one must add $n-4$ loops at the basepoint.  The darkest edges represent a choice of tree $\mathcal{T}$ lifting a maximal tree in the 1-skeleton of $K_2/Aut(F_n)$, and the lighter solid edges represent the additional edges in $\mathcal{E}$.   For $n=3$, the picture is the same except that the backmost triangle is missing and every remaining graph has one fewer loop at the basepoint.  For $n=2$, there are only the three leftmost triangles and every remaining graph has two fewer loops at the basepoint.

%\begin{figure}[ht]
%\centering
%\input{quotienttree.pstex_t}
%\caption{A lift of the quotient $K_2/Aut(F_4)$.  For $n>4$, add loops to the basepoint.}
%\label{qt}
%\end{figure}

\subsection{Vertex stabilizers}
We will use the labels in Figure~\ref{qt} to refer to the vertices of $\mathcal{T}$, and we will denote by $W_{n-i}$  the subgroup of $W_n$ which permutes and inverts the last $n-i$ generators.  

The stabilizer of a vertex  in the spine  of Auter space can be identified with the automorphism  group of the marked combinatorial graph associated to that vertex \cite{vogtmann:survey}.  We compute:
\\
\\
$G_{0} = stab(v_0) = W_n$.\\
\\
$G_{1} = stab(v_1)  = \Sigma_3 \times W_{n-2}$, where $\Sigma_3$ is the symmetric group on the three vertical edges, which is generated by $\sigma_{12}$ and $\eta$.  \\
\\
$G_{4} = stab(v_4) =  (\Z_2 \times \Z_2) \times W_{n-2}$.  Here one $\Z_2$ is generated by $\sigma_{12}$ and the other $\Z_2$ is generated by $\tau_1 \tau_2$.\\
\\
Since $v_3=\eta v_4$, 
$G_{3} =  stab(v_3)= \eta G_{4}\eta$, which is generated by $\eta\tau_1\tau_2\eta$ and $\eta\sigma_{12}\eta$.\\
\\
$G_{2} = stab(v_2) = (\Z_2 \times \Z_2)\times W_{n-2} $.  The first $\Z_2$ is generated by $\eta\tau_1\tau_2\eta$ and the second by $\tau_1$. 
\\
\\
$G_7=stab(v_7) =  D_8 \times W_{n-3}$, where $D_8$ is the dihedral group generated by 
$\eta\sigma_{12}\eta$ and $\tau_2\sigma_{13}$.  
\\
\\
$G_{6} = stab(v_6) = D_8 \times W_{n-3}$.  Here $D_8$ is the dihedral group generated by $\sigma_{12}$ and $\eta\sigma_{13}\tau_2\eta$. Note that $G_7= \eta G_6\eta$, since $v_7=\eta v_6$.
\\
\\
$G_{5} =  stab(v_5)=  \Sigma_4\times W_{n-3}$.  The symmetric group $\Sigma_4$ corresponds to permuting the edges which are not loops and is generated by the involutions $\sigma_{12}, \eta\sigma_{13}\tau_2\eta$ and $\sigma_{23}.$  Note that $G_6$ is a subgroup of $G_5$.  
\\
\\
\noindent $G_{8} = stab(v_8)= ((\Sigma_3 \times \Sigma_3)\rtimes \Z_2) \times W_{n-4}$.  The  factor $\Z_2$ is generated by $\omega = \sigma_{13} \sigma_{24}$, the first $\Sigma_3$ is equal to $G_1$ and the second $\Sigma_3$ is $\omega G_1 \omega$ .

\begin{figure}[ht!]
\labellist
\small\hair 2pt 
\pinlabel{$G_0\cong W_n$}   at 157 -10
\pinlabel{$G_1\cong\Sigma_3$} at 157 135
\pinlabel{$G_2\cong \Z_2\times\Z_2$} at 15 -10
\pinlabel{$G_3\cong \Z_2\times\Z_2$} at 15 160
\pinlabel{$G_4\cong \Z_2\times\Z_2$} at 37  225
\pinlabel{$G_5\cong \Sigma_4$} at 320 -10
\pinlabel{$G_6\cong D_8$} at 320 160
\pinlabel{$G_7\cong D_8$} at 280 225
\pinlabel{$G_8\cong (\Sigma_3\times\Sigma_3)\rtimes \Z_2$} at 171 235
\scriptsize
\pinlabel{$\langle\sigma_{12}\rangle$} at 140 90 %e_01
\pinlabel{$\langle\tau_1\rangle$} at 70 12 %e_02
\pinlabel{$\langle1\rangle$} at 70 72  %e_03
\pinlabel{$\langle\tau_1\tau_2,\sigma_{12}\rangle$} at 40 185 %e_04
\pinlabel{$\langle\sigma_{12},\sigma_{23}\rangle$} at 235 12 %e_05
\pinlabel{$\langle\sigma_{12}\rangle$} at   255 72 %e_06
\pinlabel{$\langle\sigma_{13}\tau_2\rangle$} at 280 185 %e_07
\pinlabel{$\langle\sigma_{12},\sigma_{13}\sigma_{24}\rangle$} [l] at 195 207 %e_08
\pinlabel{$\langle\eta\sigma_{12}\eta\rangle$} at  60 140 %e_13
\pinlabel{$\langle\sigma_{12}\rangle$} at 95 200 %e_14
\pinlabel{$\langle\sigma_{12}\rangle$} at 260 140 %e_16
\pinlabel{$\langle\eta\sigma_{12}\eta\rangle$} at 230 175 %e_17
\pinlabel{$\langle\eta,\sigma_{12},\sigma_{34}\rangle$} at 160 207 %e_18
\pinlabel{$\langle\eta\tau_1\tau_2\eta\rangle$} [r] at 5 72 %e_23
\pinlabel{$\langle\sigma_{12},\eta\sigma_{13}\tau_2\eta\rangle$} [l] at 315 72 %e_56
\endlabellist
\centering
\vskip 10pt
\includegraphics[scale=0.9]{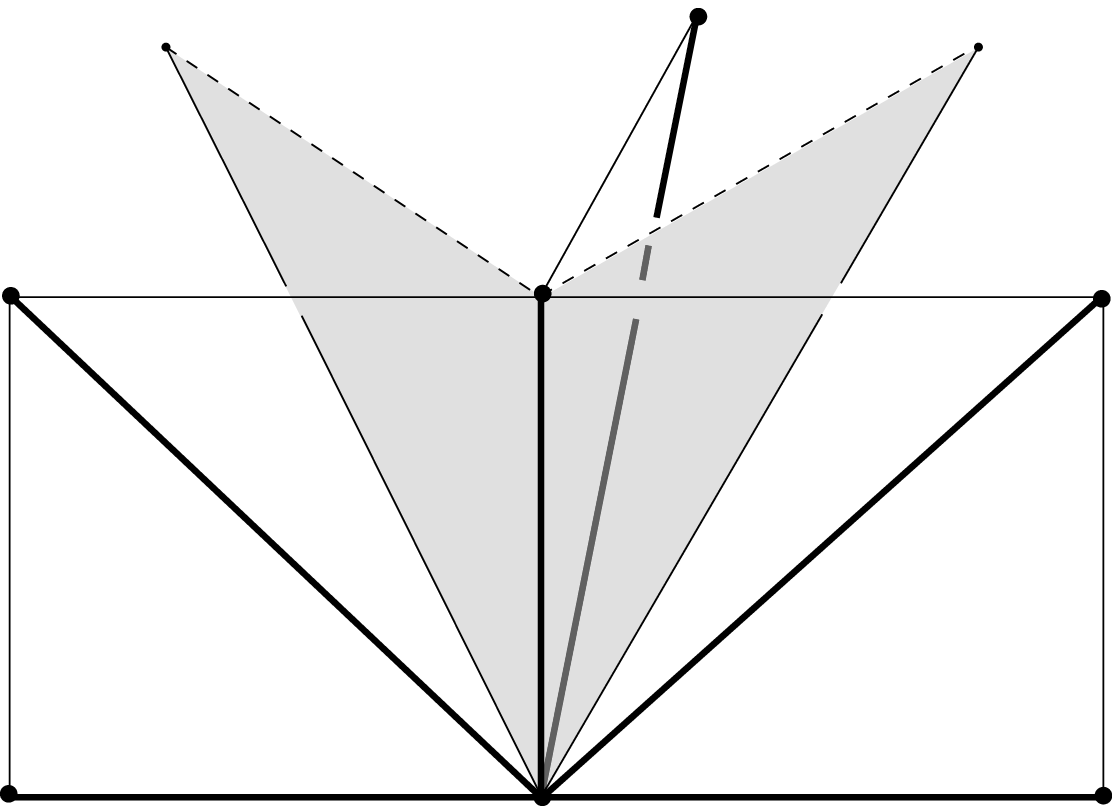}
\vskip 20pt
\caption{Edge and vertex stabilizers with $W_{n-k}$ factors omitted, except at $G_0$.  The vertex stabilizers are generated by the incoming edge stabilizers. }
\label{qtwo}
\end{figure}

By Brown's theorem, $Aut(F_n)$ is generated by the vertex stabilizers $G_i$ corresponding to vertices of $\mathcal{T}$, i.e. $G_0, G_1, G_2, G_3, G_5, G_6$ and $G_8$,  together with a generator $t_{e}$ for each of the 13 edges of $\mathcal{E}$.  We denote the oriented edge from $v_i$ to $v_j$ by $e_{ij}$.  

\subsection{Tree relations} If $e\in \mathcal{T}$ (i.e. $e=e_{0k}$ for $k\in\{3,6,8\}$ or $e_{k0}$ for $k\in\{1,2,5\}$), then the tree relations set $t_e=1$.  

\subsection{Face relations} If all edges of a triangle $\Delta$ are in  $\mathcal{E}$  and two of the edges lie in $\mathcal{T}$, then $h=1$ and $g_e=1$ for all edges in the boundary of $\Delta$, so
the face relation associated to $\Delta$ reduces to $t_e=1$ for the third edge $e$.  We now have $t_e=1$ for all edges except $e_{04}$ and $e_{07}$.  

 The only faces which do not have two edges in $\mathcal{T}$ are the shaded faces labeled $\Delta_{104}$ and $\Delta_{107}$ in Figure~\ref{qt}.  
The boundary of $\Delta_{107}$ is given by the edge-path loop $e_{10}e_{07}e_{71}$.
The first edge $e_{10}$ is in 
$\mathcal{T}$,  giving $h_{10}=1$ and $g_{e_{10}}=1$, 
so $g_{10}=1$.  The second edge $e_{07}$ is in $\mathcal{E}$, so $h_{07}=1$;  this edge has $t(e_{07})=v_7=\eta v_6$,  so $w(e_{07})=v_6$ and $g_{e_{07}}=\eta$, giving $g_{07}=\eta$.  The last edge  $e_{71}$ is equal to $\eta e_{61}$, and $e_{61}\in\mathcal{E}$.  Thus $h_{71}=1$, and $g_{e_{71}}=1$, giving $g_{71}=1$. We have $g_{10}g_{07}g_{71}=\eta$, and the relation associated to ${\Delta_{107}}$ is now
$$1\cdot t_{e_{10}}\cdot 1\cdot t_{e_{07}}\cdot 1\cdot t_{e_{71}}=\eta$$
which reduces to $t_{e_{07}}=\eta$.  An identical calculation for $\Delta_{104}$ gives $t_{e_{04}}=\eta$, since $\eta v_4=v_3$. 

\subsection{Edge relations}  The edge relations identify generators of the vertex groups with the appropriate products of the $\tau_i$, $\sigma_{ij}$ and $\eta$.  
If $e\in \mathcal{T}$,  the edge relations identify all of the generators written as products of $\tau_i$ and $\sigma_{ij}$ in our descriptions of the $G_i$  with the corresponding elements of $G_0=W_n$; in particular, the subgroups $W_{n-k}$ of the stabilizers are all identified with the corresponding subgroup  of $G_0$.   
The edge relation associated to $e_{04}$ identifies the generators of $G_3$ with $\eta\tau_1\tau_2\eta$ and $\eta\sigma_{12}\eta$, since $\eta v_4=v_3$.   The edge relation associated to $e_{23}$ then identifies the first generator of $G_2$ with $\eta\tau_1\tau_2\eta$.

The edge relation associated to $e_{07}$ identifies the second generator of $G_6$ with $\eta\tau_2\sigma_{13}\eta$, since $\eta v_7=v_6$.  The edge relation associated to $e_{56}$ identifies $G_6$ with the subgroup of $G_5$ generated by $\eta\tau_2\sigma_{13}\eta$ and $\sigma_{12}$.  

The edge relation associated to $e_{18}$ identifies the first $\Sigma_3$ with the corresponding $\Sigma_3$ subgroup of $G_1$.

\subsection{Stabilizer relations}  We will not list the relations in $G_0=W_n$.
The relations in $G_1$ which do not come from $G_0$ are those involving $\eta$, i.e.
 \begin{eqnarray}
\eta^2=1\\
 (\sigma_{12}\eta)^3=1 \\
 (\eta \tau_{i})^2= 1  & i>2\\
(\eta \sigma_{ij})^2= 1  & i,j>2.
\end{eqnarray}
Since $G_3=\eta G_4\eta$, and $G_4$ is a subgroup of $W_{n}$, $G_3$ does not contribute any new relations. The fact that the generators of $\Z_2\times \Z_2\leq G_2$ commute contributes the relation $(\eta\tau_1\tau_2\eta\tau_1)^2=1$, which looks a little nicer if we conjugate by $\eta\tau_1$:
 \begin{eqnarray}
 ((\eta\tau_1)^2\tau_2)^2=1.
\end{eqnarray}
In the dihedral group $D_8\leq G_6$, the fact that the product of our generators has order 4 contributes a new relation
 \begin{eqnarray}
 (\eta\sigma_{13}\tau_2\eta\sigma_{12})^4=1.
\end{eqnarray}
%It is straightforward to check that this relation plus the relations inherited from $G_0$ determine  $G_5$, so $G_5$ contributes no new relations.  NO IT's NOT!!!
The symmetric group $\Sigma_4 \leq G_5$ is generated by the involutions $\sigma_{12},\sigma_{23}$ and $\phi=\eta\sigma_{13}\tau_2\eta$, with relations $(\sigma_{12}\sigma_{13})^3=1, (\sigma_{12}\phi)^4= 1$ and finally $\sigma_{12}\phi\sigma_{12}(\sigma_{23}\phi)^2=1$. The first relation comes from $G_0$ and the second from $G_6$,  so $G_5$ adds only the third relation, i.e.
 \begin{eqnarray}
 \sigma_{12}\eta\sigma_{13}\tau_2\eta\sigma_{12}(\sigma_{23}\eta\sigma_{13}\tau_2\eta)^2= 1.
 \end{eqnarray}

The fact that the two copies of $\Sigma_3$ which are contained in  $G_8$ commute produces the relation $(\sigma_{14}\sigma_{23}\eta\sigma_{23}\sigma_{14}\eta)^2=1$, i.e. 
\begin{eqnarray}
 (\sigma_{14}\sigma_{23}\eta)^4=1.
\end{eqnarray}
All other relations in $G_8$ are consequences of this and relations in $G_0$ and $G_1$; for example the fact that $\eta\sigma_{12}\eta=\sigma_{34}$ commutes with $\eta$ is a consequence of relations already accounted for in $G_1$.     

This completes the proof of Theorem~\ref{presentation} and its corollaries.

\medskip\noindent{\bf Acknowledgments:}  
Karen Vogtmann was partially supported by NSF grant DMS-0204185.

\end{document}